\documentclass[10pt]{amsart}
\usepackage{latexsym}
\usepackage{amssymb}
\usepackage{amsmath}
\usepackage{endnotes}

\begin{document}

\title{The Tragedy of Mathematics in Russia}
\author{S.~S. Kutateladze}

\address[]{
Sobolev Institute of Mathematics\newline
\indent 4 Koptyug Avenue\newline
\indent Novosibirsk, 630090
\indent RUSSIA}
\email{
sskut@math.nsc.ru
}
\begin{abstract}
This is a brief overview of the role of mathematicians in the so-called
``Luzin Case''
as well as some analysis of the mathematical and humanitarian roots of the affair.
\end{abstract}
\date{March 8, 2007---May 9, 2011}

\maketitle



Nikola\u\i{} Nikolaevich  Luzin (1883--1950) was one of the founding fathers of
the Moscow mathematical school. The list of his students
contains Full Members of the Academy P.~S.
Aleksandroff (1896--1982), A.~N. Kolmogorov (1903--1987),
M.~A. Lavrentiev (1900--1980), P.~S. Novikov (1901--1975);
Corresponding Members  L.~A. Lyusternik (1899--1981),
A.~A. Lyapunov (1911--1973), D.~E. Menshov (1892--1988),
A.~Ya. Khinchin (1894--1959), L.~G. Shnirelman  (1905--1938);
and many other mathematicians.

The outstanding roles in the development of mathematics in Siberia
were performed by Lavrentiev and Lyapunov, Luzin's direct descendants
as well as Academicians A.~I. Maltsev (1909--1967) and A.~A. Borovkov,
students of Kolmogorov.

\section*{The Case  Against Luzin}

Tsunami swept over the Russian mathematical community
in 1999 after publication of the complete shorthand notes
of the meetings of the notorious emergency Commission of
the Academy of Sciences of the USSR on the
case of Academician Luzin~\cite{Case}. Soon the article \cite{Lorentz} appeared in
the USA which  revealed  the personal
testimony of
G.~G. Lorentz (1910--2006)
about the mathematical life of that time in the USSR.\endnote{I am
very grateful to  Professor W.~A.~J. Luxemburg for
attracting my attention to the inadvertent omission
of a reference to \cite{Lorentz} in the draft
of this paper.}

The Commission for the ``hearing of the case of Ac[ademician] Luzin'' was convened by the Presidium of the Academy of Sciences
of the USSR after the article ``Enemies under the  Mask of  a Soviet Citizen''
in the {\it Pravda} newspaper on July~3, 1936.

Luzin was  accused of all theoretically possible instances of
misconduct in science and depicted as an enemy that combined
``moral unscrupulousness and scientific dishonesty
with deeply concealed enmity and hatred to every bit of the
Soviet life.'' It was alleged that he publishes ``would-be scientific papers,''
``feels no shame in declaring the discoveries of his students to be
his own achievements,'' and stands close to the ideology of the ``black hundred'',
orthodoxy, and monarchy ``fascist-type modernized but slightly.''
The closing  of the lampoon read:

\begin{itemize}
\item[]{\small
The Soviet scientific community tears away
from you the mask of an honest scientist, leaving you in the altogether,
and so you appear before the eyes of the world as a paltry individual
who pretends to champion ``pure science'' but betrays the interests of science,
merchandizing it to appease you former bosses---the present-day masters,
of faschistoid science. The Soviet community will perceive the story
about Academician Luzin as another object-lesson of the fact that
the adversary never lays down his arms, that he camouflages himself
more  skilfully, that the methods of his mimicry becomes more
diverse, and that vigilance remains the most demanded trait of every
Bolshevik and every Soviet citizen.
}
\end{itemize}

All Russian scientists of the elder generation knew about the {\it Pravda\/}
editorial and the savage dissolution of ``Luzinism.''
There was no denying that the initiation of the campaign for
discrediting Luzin was carried out by the symbiosis of
the party and repressive machinery of the USSR.
Behind the scenes of the campaign loomed the grim figures
of E.~N.  Kolman  (1892--1979)
and  L.~Z. Mekhlis (1889--1953),  typical representatives  of the {\it Oprichnina\/}
of the Stalinist epoch.
The  former was Head of the Department of Science of
the Moscow Committee of the All-Union Communist Party (Bolsheviks), while
the latter was Editor-in-Chief of  {\it Pravda}.

The Luzin case had been considered only within the
context of the general crimes of the Stalin totalitarianism for a long while.
The newly-published archive files have opened to the public the new
circumstances---some students of Luzin were rather active participants of
the political assault on  their teacher.
The key role among them was played by P.~S. Aleksandroff
who headed  the Moscow topological school.

S.~P. Novikov\endnote{Serge\u\i{} Petrovich Novikov  is a Fields medalist,
the son of L.~V. Keldysh and P.~S. Novikov who were  students of Luzin.}
 wrote in~his reminiscences (cp.~\cite[p.~45]{Novikov}):

\begin{itemize}
\item[]{\small
Investigation was launched  that time by
my father (seemingly, with Lyusternik and Lavrentiev, the latter knowledgable
about the Party circles).
They found out  there was a letter from
P.~Aleksandroff to some influential guy,
Khvorostin by name.\endnote{G.~K. Khvorostin (1900--1938) was an educated
mathematician, Rector of Saratov State University in 1935--1937.
A.~Ya. Khinchin  qualified him as ``the Party youth''
(cp.~\cite[p.~100]{Case}).}
This letter narrated the nasty  deeds of Luzin. Khvorostin
resided in Saratov and had great connections within the C[entral] C[ommittee of the Party].
He hated Luzin, which was known to everyone. It was Khvorostin as they
guessed out who submitted materials  to the CC and initiated the article.
Pavel Sergeevich was a great master of a billiard shot!}
\end{itemize}

Luzin was especially annoyed by the invectives of P.~S. Aleksandroff
which were aimed at the declension of Luzin's contribution to
analytic set theory (cp.~\cite{Lorentz_Sets}).
It is now in common parlance to call an analytic set the continuous image of a Borel
subset of the reals. These sets are often associated with the names  Aleksandroff and Suslin
and called $A$-sets or Suslin sets.
Note that P.~S. Aleksandroff commented at the meeting of the Commission on July 6, 1936
as follows (cp.~\cite[p.~90]{Case}):

\begin{itemize}
\item[]{\small
Suslin called them $A$-sets. But he never told that he called them so in my honor.}
\end{itemize}

In his reminiscences of 1979
he claimed quite the opposite (cp.~\cite[p.~235]{PSA}):

\begin{itemize}
\item[]{\small
It was exactly then that Suslin proposed to call the
new set theoretic operation that I had constructed the A-operation while calling $A$-sets those that
result from application of the operation to closed sets.
Doing so, he emphasized  that he was suggesting this terminology
by analogy with Borel sets which were customarily called $B$-sets those days.}
\end{itemize}

There was a mysterious relevant episode about the history
of the discovery of $A$-sets which was narrated by A.~P. Yushkevich (1906--1993)
 not long before his death (cp.~\cite{Yush_meetings}):
\begin{itemize}
\item[]{\small During our rather often meetings in the second
half of the 1970s, I  repeatedly asked P.~S. about the history of
the discovery of $A$-sets. One day  he rang me up and suggested that
I would come to his apartments  at an appointed time when he would stay alone.
I had to enter into the waiting room without ringing the door bell, and P.~S.
would be seated on a chair quite opposite to the door. All happened
exactly as planned. His narration was lengthy. I wrote it down.
On the next day I brought the typewritten text of my record, and P.~S. told that
it conforms completely with what he had, but stipulated the condition that
the record be never published during his life time.
When P.~S. had passed away in 1982,
I decided to take the advice of A.~N. [Kolmogorov ({\bf S.~K.})].
After  listening the text, A.~N. told me that to publish it  was still too early.
Now there are neither P.~S. nor A.~N., and I have nobody else
to take advice of. I decided to seal the narration of P.~S.
Aleksandroff in an envelope and put it to the Archive of
the Academy of Sciences of the USSR for deposition in the Aleksandroff Fond
with the inscription: ``Received from A.~P. Yushkevich.
Open only after  10 years since his death.''\endnote{
 J.-M. Kantor informed me that his attempts, and the attempts by
S.~S. Demidov either, were in vain to find this envelope in the Archive of the
Russian Academy of Sciences ({\bf S.~K.}).}}
\end{itemize}

Rather active at the meetings of the Commission were A.~N. Kolmogorov,
L.~A. Lyusternik, A.~Ya. Khinchin, and L.~G. Shnirelman. The political attacks on
Luzin were vigorously supported by members of the Commission
S.~L. Sobolev (1907--1989) and O.~Yu. Schmidt (1891--1956).
A.~N. Krylov (1863--1945) and S.~N. Bernstein (1880--1968) revealed valor in
the vigorous defence of Luzin.
The  final clause of the official Resolution of the Commission
read as follows (cp.~\cite[p.~296]{Case}):

\begin{itemize}
\item[]{\small
Everything of the above, summarizing the overwhelming
material evidence in possession of the Academy of Sciences,
completely ascertains the characteristics of Luzin in the {\it Pravda\/}
newspaper.}
\end{itemize}

All participants of the events of 1936 we discuss had left this world.
They seemingly failed to know that the files of the Commission
are all safe and intact. Today we are aware in precise detail
of what happened at the meetings of the Commission and around the whole case.
The mathematical community painfully reconsiders the events and rethinks
the role of the students of Luzin in his political execution.

\section*{Roles of Luzin's Students}

P.~S. Novikov and M.~A. Lavrentiev  were not listed
as participants of the public persecution of Luzin
(despite the fact that both were  mentioned at the meetings
of the Commission among the persons robbed by Luzin).
It transpires now why M.~A. Lavrentiev was the sole author of
a~memorial article~\cite{MAL} in  {\it Russian Mathematical Surveys\/}
on the occasion of the ninetieth anniversary of the birth of Luzin.
He also included this article in the collection \cite{MAL80} of his papers on
the general issues of science and life. M.~A.~Lavrentiev was the chairman of the
editorial board of the selected works of Luzin which were
published by the decision of the Academy of Sciences of the USSR
after the death of Luzin on the occasion of the 70th anniversary of
his birth. P.~S. Aleksandroff and A.~N. Kolmogorov were absent from
the editorial board.

Practically the same are the comments on their relationship with Luzin  which were left
by  P.~S. Aleksandroff and A.~N. Kolmogorov. Their statements are
still shared to some extent by their numerous students.
It is customary to emphasize that Luzin was not so great a mathematician
as his students that had persecuted him. Some moral fault is persistently incriminated
to Luzin in the untimely death of M.~Ya. Suslin (1894--1919) from  typhus fever.
Luzin is often blamed for all his disasters at least partly.
He is ascribed such traits of character as theatricality,
envy of others' success, hypocrisy, plagiarism, and inclination to
intriguing.\endnote{No person with these defects of personality
could ever become the founder of ``Lusitania''---the most successful mathematical
school in the history of world science. Therefore, there exists as rather likely theory
of ``two Luzins'': one of the epoch of Lusitania and the other of the times of the Luzin case.}
He is said to deserve all  punishments and if not all then
it is not his students' fault but stalinism and the curse of the epoch.
These arguments reside in the minds of not only  the elders but also  the
youngsters.  The best of them view the Luzin case as the mutual
tragedy of all participants.

We should however distinguish the personal tragedy of Luzin from
the tragedy of the Moscow school and the tragedy of the national mathematical community.
The students of  Luzin who participated in the persecution of the teacher
never considered  their own fates tragical.

P.~S. Aleksandroff wrote in his reminiscences (cp.~\cite[p.~90]{PSA}):

\begin{itemize}
\item[]{\small
Knowing  Luzin in his green creative years, I got acquaintance
with a truly inspired teacher and scholar who lived
only by science and in the name of science.
I met a person who resided in the sphere of the sublime human treasures
which is forbidden for any rotten ghost or spirit.
When a human being leaves this sphere (and Luzin had left it once),
he is doomed to surrender to the forces that were described by Goethe
as follows:}

\qquad\vbox{
\item[]{\small
Ihr f\"uhrt in's Leben uns hinein,\\
Ihr lasst den Armen Schuldig werden\\
Dann \"uberlasst Ihr ihn der Pein,\\
Denn jede Schuld r\"acht sich auf Erden.\\
Into our life you lead us in,\\
The wretch's guilt you bring to birth,\\
Then bring affliction down on sin,\\
For all guilt takes revenge on
Earth.\endnote{
P.~S. Aleksandroff cited the poem {\it Harfenspieler} dated as of 1795
by Johann Wolfgang von Goethe (1749--1832) and gave
a~rough translation into Russian. The lines in English  here belong
to Vernon Watkins (1906--1967).}
}
}
\item[]{\small
\noindent
In his terminal years Luzin saw the bottom of the sour bowl of
the revenge that was described by
Goethe.}
\end{itemize}

It is worth observing that A.~Ya. Khinchin, hostile to Luzin,
commented on the accusations that Luzin drove
M.~Ya. Suslin to death (cp.~\cite[p.~391]{Yush_meetings}):

\begin{itemize}
\item[]{\small
Suslin is called the student perished
by N.~N. Luzin. Why, when a man dies from typhus fever
this is a rather exaggerated expression.
In fact Suslin could  possibly get typhus fever in Ivanovo.
Furthermore, in the common opinion it was N.~N.
who tried and expelled Suslin from Ivanovo.
But the transfer from Moscow to Ivanovo
I view as a favor to Suslin who was not hostile to
Luzin in those days.}
\end{itemize}

\noindent
Narrating his reminiscences of  P.~S. Aleksandroff,  A.~N. Kolmogorov
told in~1982 (cp.~\cite[p.~10]{Tikh}):
\begin{itemize}
\item[]{\small
My entire life as a whole was full of happiness.}
\end{itemize}

\noindent
 Neither he nor
 Aleksandroff nor other participants of the persecution of Luzin
had ever treated the
Luzin case as a common tragedy with Luzin. They were correct in this judgement
 but on the grounds completely different from those they declared.

If Luzin were guilty then his fault would belong to the sphere of
the personal mathematical relations between a teacher and a student.
No convincing evidence of Luzin's plagiarism was ever
submitted. The alleged accusations that he ascribed
to H.~Lebesgue (1875--1941) or kept a grip of
Suslin's  results  are poorly disguised and baseless.
To prove the scientific misconduct of Luzin it was alleged that
Luzin ingratiated and flattered H.~Lebesgue by ascribing
Luzin's sieve method to H.~Lebesgue.
On the other hand, H.~Lebesgue wrote in his preface to the Luzin book
on analytic sets (cp.~\cite{Lebesgue}):

\begin{itemize}
\item[]{\small
Anyone will be astonished
to find out from Luzin's book that I had incidentally invented
the sieve method and was the first to construct an analytic set.
But nobody could be more amazed than me. Mr. Luzin
feels himself happy only when he has managed to ascribe his own discoveries to
someone else.}
\end{itemize}

\noindent
The students were ``more Catholic than
the Pope.''\endnote{
The corresponding excerpt of the shorthand minutes of the meeting of the Commission
on July 13, 1936 is as follows:
\cite[p.~196--197]{Case}:
\hfill\break
\tiny
\indent
{\bf Aleksandroff.} As regards obsequiousness,
I suggest that we  will use the genuine words of the lips of
Lebesgue: (Reads in French). Concerning this matter, I have
explanations that I can explicate as thoroughly as need be.
The  ``strange mania'' in question, I would say, is a deeply
premeditated idea. He ascribed to Lebesgue his belongings
and he did it in so dopey manner. No sane person would ever
ascribe them to Lebesgue. The thing is that doing so
he creates his reputation of the person who
ascribes his own ideas to someone else; but when
the matter concerns his own students then he robbed their
belongings while hiding behind this screen.
\hfill\break
\indent
{\bf Lyusternik.} This kind of defence sounded at
the meeting in the [Steklov {\bf(S.~K.)}] Institute.
Precisely this way of defence that was clearly inspired by him:
How might it happen that N.~N. grips the results of the others
if Lebesgue himself writes these words about him?
\hfill\break
\indent
{\bf Aleksandroff.} This is an obsequious system
since it is uncustomary in academic circles to ascribe  someone's own
results to anybody else. Therefore, we see here, on the one hand,
his flattery of Lebesgue and, on the other hand,
the arranging of the screen that allows him to behave so.
}

In any case, we cannot help but ascertain that there clearly was a conflict of
generations---the precipice of alienation and misunderstanding between Luzin and his most
successful students.
It is easy to assume the genuine or imaginary injustice and prejudice
of Luzin in citing his students as well as the genuine or
imaginary feebleness of Luzin in overcoming mathematical obstacles.
We may agree to see hypocrisy in Luzin's decision to vote against
P.~S. Aleksandroff in the elections to a vacancy of an academician
despite his personal letter of support of P.~S. Aleksandroff  to  A.~N. Kolmogorov.
Well, there is nothing untypical of the academic manners or extraordinary
in Luzin's  conduct, is there?
It is the true background of the Luzin case, isn't it?

Available is the following testimony of W.~Sierpi\'nski
(1882--1969), a famous Polish mathematician who was declared to be
a  ``blatant black hundredist'' at the meetings of the Commission of
the Academy of Sciences of the USSR on
the Luzin case (cp.~\cite[p.~124]{Dugac}):

\begin{itemize}
\item[]{\small
In his letter as of June 27, 1935---which was a year
ago---Mr.~Luzin wrote:
``Returning now to my very difficult self-defence
against  the ascription to Suslin of the results which
he had no rights to   and which were absent even in
his thoughts, I must say that this self-defence was provoked by
a very grave and absolutely impending danger. Mr.~Aleksandroff
has dreams of entering the Academy of Sciences as a full member
by dismissing me. To this end he requests  that my
contributions be reconsidered, claiming that I has no right
to be a member of the Academy since all my ideas are stolen from Suslin.
This reconsideration is rather likely and feasible.''
When I was in Moscow in September, 1935,
Mr. Aleksandroff assured me that the apprehensions of Luzin
are purely imaginary and that he respects Luzin, his former teacher.
In my presence  Aleksandroff shook hands with Luzin and declared that he
would always be a friend of Luzin.}
\end{itemize}

The pretentious reconciliation of  P.~S. Aleksandroff with~Luzin
which was described by W.~Sierpi\'nski and which was later publicly
refuted by  P.~S. Aleksandroff is in no way similar to the refusal of Luzin
to support the election of P.~S. Aleksandroff as an academician, isn't it?
It is in general belief that this refusal was the reason for
A.~N. Kolmogorov  to  slap the face of Luzin publicly
in~1946.\endnote{About various versions of this episode
see the recent book by L.~Graham and J.-M. Kantor \cite[p.~186]{Naming} as well as
the reminiscences of S.~M. Nikolski\u\i{} \cite[p.~155]{Nikolsky} and
S.~P. Novikov \cite[p.~22]{Novikov}.
}
L.~S. Pontryagin (1908--1988)
wrote on December~24, 1946 (cp.~\cite[Letter No.~49, pp.~89--91]{Gordon}):

\begin{itemize}
\item[]{\small
You are asking about cooperation of Kolmogorov and Luzin.
Surely, this needs to be narrated rather that written since the mode of voice is essential to
render all properly. In summer Kolmogorov told me that his only concern
about the election of Aleksandroff consists in the fact that
he became an obvious candidate  four months before the elections.
{\it Pusiks}\endnote{With the letter {\tt u} pronounced as {\tt oo} in {\tt soon}, {\it Pusiks\/}
was the collective equivocal
nickname of P.~S. Aleksandroff and A.~N. Kolmogorov. The singular {\it Pusik\/} was applied only
to P.~S. Aleksandroff. For the  English ear the word {\it Pusiks}
has an obscene connotation, but this effect is completely absent in Russian.}
carried out a great preparatory work in the sense that
entered into various agreements with academicians. For instance, there was a promise to
Vinogradov that Kolmogorov would support Lavrentiev in the event that Vinogradov
would support Aleksandroff. All in all  it seemed that everyone will vote for
Aleksandroff. For example, Bernstein had nominated Aleksandroff at a meeting of the
[Steklov] Institute as well as Chebotar\"ev, I must mention for then sake of truth.
Kolmogorov had made an agreement with the bosses of the Academy that he will be
a member of the Board of Experts.\endnote{This is an interim committee convened at the elections
that discusses all nominees and chooses those that are recommended to be elected.
The names of the recommended nominees appear at the top on  voting bulletins, which
prompts the choice of those unacquainted with the candidates.} The first annoying news was the fact
that he had not been assigned to the Board, but he still hoped that this was not essential already.
After the session  of the Board of Experts there were a few
closed meetings of all academicians where they discussed the candidates, and only at that stage
Kolmogorov became aware that none of the members of the Board of Experts had supported
Aleksandroff. On the contrary, Bernstein vigorously objected against him,
claiming that Aleksandroff's area of research is harmful. Bernstein's behavior
still seems highly illogical for me; probably he had just quarreled with
{\it Pusiks}. All the rest is rather comprehensible.
Lavrentiev turned out an unquestionable candidate somehow and he had no need in the
support of Kolmogorov who was not a member of the Board at that.
Therefore, Vinogradov had no need in Kolmogorov and  {\it Pusik}.
As regards Sobolev and Khristianovich, the former deeply hated {\it Pusik}
for expelling Sobolev from the directorship [of the Steklov Institute ({\bf S.K.})],
while  the latter is a crony and companion of Sobolev.
In these circumstances there was practically  no hope of success.
The only remaining possibility was that some  academicians among mathematicians
would support Aleksandroff, for physicists wanted to support him but surely
that could not confront all mathematicians.
Luzin became the  hope of
{\it Pusiks}
He was invited to Komarovka
and promised his support.
However, he spoke against Aleksandroff at the final closed meeting.
Departing from this meeting,  Kolmogorov was  absolutely upset and stung.
He came to Luzin and said that he would have nothing in common
with Luzin ever since. Luzin pretended that he did not
understand anything and began to talk as follows:
``Dear me, calm down. Forget it. You are ill. Relax.''
This is what must be narrated with expression.
Kolmogorov then answered him: ``So what shall I do to you:
spit at your physiognomy or slap your mug?''  After a~short thought, he dared
the latter.\endnote{
The elections to the Academy of Sciences of the USSR in 1946 took place
on November 30. It was M.~A.~Lavrentiev and  I.~G.~Petrovski\u\i{} (1901--1973) who were
elected to fill the  mathematical vacancies of a full member  in the Division of Physics
and Mathematics.
A.~D. Alexandrov (1912--1999), N.~N. Bogolyubov (1909--1992),
L.~A. Lyusternik, and V.~V. Stepanov (1889--1950)
became corresponding members of the Academy.
}}
\end{itemize}

V.~M. Tikhomirov disclosed  the correspondence between Luzin and A.~N. Kolmogorov on the eve
of the elections of  P.~S. Aleksandroff.
In the fall  of 1945 Luzin wrote to A.~N. Kolmogorov (cp.~\cite[p.~80]{Tikh}):

\begin{itemize}
\item[]{\small
Now about another matter: the time is coming of the elections to the Academy.
It would be an utmost injustice, had these be happened without Pavel Sergeevich.
His works  whose echoes are met  throughout the world literature,
his splendid ripe years---the completion of maturity and wisdom and he himself as
the man of perfect attraction---all these made us see in him a wordy candidate
who activity is invaluable for the Academy.}
\end{itemize}

A.~N. Kolmogorov overestimated the position of Luzin
who had written nothing more  than that he will support
P.~S. Aleksandroff as a nominee for the elections (which Luzin  fulfilled in due time).

In the reply letter of October 7, 1945  A.~N. Kolmogorov
remarked (cp.~\cite[p.~82]{Tikh}).

\begin{itemize}
  \item[]{\small
 Since for many years I am engrossed in  preventing  the various random and immaterial
 circumstances to hinder the election of Pavel Sergeevich which would be
 a rather just event in my opinion; therefore, I indeed appraise very highly your
 readiness to support all necessary actions for the success when this is really desired.
 }
\end{itemize}

 Clearly, A.~N. Kolmogorov decided without proper  grounds that
 Luzin would support the election rather than  nomination of P.~S. Aleksandroff.
 Luzin however promised nothing of the sort in his letter. It is the tradition of a long
 standing that nomination and election to the Academy of Sciences and other similar
 institutions are sufficiently independent procedures.

 Usually A.~N. Kolmogorov was viewed as a calm person not
 liable to fits and  extremes of temper. Therefore, he seemingly needed
 some special provocation from Luzin for slapping in Luzin's face,
 which led to some apocrypha about an obscene remark from Luzin at the
 elections of~1946.\endnote{E.g., see  \cite[p.~186]{Naming} and  \cite[p.~22]{Novikov}.}
An analogous version was mentioned by V.~I. Arnold in private correspondence.
It is not excluded that the available hints on {\it topolozhestvo}\endnote{This is a conjunction
of the Russian words for topology and pederasty---something like ``topologasty'' in English.}
is a produce of the 1950s  put in gossips for rehabilitation of the instigators of
the ``Luzin case.''
 But the students of A.~N. Kolmogorov indicated one quite unknown trait of
his personality.  V.~M. Tikhomirov wrote (cp.~\cite[p.~83]{Tikh}):

\begin{itemize}
  \item[]{\small
  \dots it should be noticed there was some rather unsavory particularity:
  sometimes he lost his temper.}
\end{itemize}

 V.~I. Arnold witnessed (cp~\cite[p.~50]{Arnold}):

\begin{itemize}
  \item[]{\small
  Andre\u\i{} Nikolaevich was never too good-natured  and he narrated not without pride about his
  clash with police at the Yaroslavl Railway Station.}
\end{itemize}

Luzin was twenty years older than  A.~N. Kolmogorov. Luzin was a teacher of
A.~N. Kolmogorov and carried the heavy burden of political accusations that were
imposed on Luzin with participation of
P.~S. Aleksandroff and  A.~N. Kolmogorov.  Luzin was granted ``mercy''
and
accepted at the country house of  A.~N. Kolmogorov and  P.~S. Aleksandroff
in Komarovka before the
elections.\endnote{V.~M. Tikhomirov wrote about the meeting in Komarovka:
``The correspondence of  L.~S. Pontryagin and his student
and friend
I.~I. Gordon reveals that Luzin was accepted and served a meal
in~Komarovka.'' Cp.~\cite[p.~83]{Tikh}.}
Everyone at the meeting remembered the most important matter that
Luzin was victimized and must surrender to the noble victors, didn't he?
It transpires now, doesn't it?
We can compare the  internal academic  matters, say Luzin's misconduct
and even plagiarism, with the accusations of  subversive  activities
against the Soviet life, can't we?
These grave and vexed questions...

Closing the meeting on July~13, 1936, Academician
G.~M. Krzhyzhanovski\u\i{} (1872--1959), Chairman of the Commission, told
in particular (cp.~\cite[p.~196]{Case}):

\begin{itemize}
\item[]{\small
We then must think over the following matter.
This fall the elections are in order and we are hinted that
there will be vacancies for 30 new academicians and 60 new
corresponding members. We are to refresh the body of the Academy,
and by the Assembly of the Academy in September
you have to ponder  over and decide the following matter:
Whom  you will recommend to be elected as corresponding members and
academicians. This will be the best outcome of the
work of this Commission.}
\end{itemize}

 No elections to the Academy were
arranged in~1936. The great elections took place
only on January 29, 1939 (cp.~\cite[No.~241, No.~242]{AN}).
The following mathematicians were elected to the
Department of Mathematical and Natural Sciences of the Academy:
A.~N. Kolmogorov and S.~L. Sobolev became full members, while
A.~O. Gelfond,  L.~S. Pontryagin, and  A.~Ya. Khinchin became corresponding members.

\section*{
Reactions of Luzin's  Contemporaries
}

All moral accusations against Luzin are rather inconvincible.
That which was submitted as proofs was inadequate even in the
times of the Commission neither for
P.~L. Kapitsa (1894--1984), nor V.~I. Vernadsky (1863--1945), nor A.~Denjoy (1884--1974),
nor Lebesgue, nor many other elder persons.

The objections of Kapitsa were expressed on July~6 in his letter to V.~M. Molotov
who was the Chairmen of the Council of the People's Commissars of the USSR.
V.~I. Vernadsky noticed in his diary on
the next day (cp.~\cite[p.~92]{VIV}):

\begin{itemize}
\item[]{\small
Letters to Luzin, Chaplygin, and Fersman about him. Majority treats
as demonstrated the slander and insinuations.  M[ay] b[e], he [is needed]
abroad but not at home. I am afraid that this disgusting article
will affect him much. Many conversations and many impressions.}
\end{itemize}

On the same day he sent a letter to Academician A.~E. Fersman (1883--1945),
a member of the Commission. V.~I. Vernadsky remarked (cp.~\cite[p.~94]{VIV}):

\begin{itemize}
\item[]{\small
I think that
such an episode would eventually be perilous to the Academy were it  led to the expulsion of
N.~N. [Luzin] from the Academy or any similar actions.
We would slide down the slippery slope.}
\end{itemize}

And on July 11  S.~A. Chaplygin (1869--1942) wrote to V.~I.
Vernadsky (cp.~\cite[p.~106--107]{Antipenko}):

\begin{itemize}
\item[]{\small
The article about Luzin is completely outrageous: Supposing that he committed a sin
of misjudging some applicant for a scientific degree or title, but how
is it possible to jump to the conclusion of sabotage from that?! \dots
As far as the accusations, slipping through the article, of fascism
and his enlisting the  old reactionary Moscow school of mathematics,
I am complete unable to understand these.
There remains the critical evaluation of Luzin's contributions. But in this regard
I must say only that this discloses the complete incompetence of the authors which
proves their minor and superficial acquaintance with his works and their
deliberate distortion of correct evaluation.
His authority is incomparable with that of Khinchin
who is counterpoised to him. But what should be done right away? How can we help N.~N.?
I did only one thing yet: I sent N.~N. a cable whose copy I attach:
``Dumbfounded by absolutely undeserved newspaper attacks against you.
Your high world-wide acknowledged  scientific authority cannot be
shaken. I hope definitely that you will find
the inner forces to face this inauthoritative criticism of your contributions calmly.
I avoid mentioning the completely groundless accusations of the other sort.}
\end{itemize}

Lebesgue's letter of August 5, 1936 is in order now.  I remind that
H.~Lebesgue was elected in 1929 to the Academy of Sciences of the USSR
for his outstanding contribution to mathematics.
The great Lebesgue, the author of that very ``Lebesgue integral'' which is indispensable
in modern mathematics,  was in the state of utmost indignation and anger.
He wrote (cp.~\cite[p.~127]{Dugac}):

\begin{itemize}
\item[]{\small
You will see that it was  not yesterday when the attacks on Luzin began with the aim
of  firing him and emptying place for Aleksandroff. You will see there
that I was already mixed in this by contrasting ``my''
science, which is bourgeoise and useless, to {\it analysis situs} [topology],
a proletarian and useful science. Since the former was the science of Luzin;
whereas the latter, the science of Aleksandroff.  What is curious is that
he begins   as Urysohn whose papers he inherited
at the same starting point that was mine. With the only difference that
Urysohn cited me  whereas  Aleksandroff has never cited me anymore
 since he must now speak badly of me in his struggle against Luzin!}
 \end{itemize}

 Another evidence of W.~Sierpi\'nski  is as follows (cp.~\cite[p.~125]{Dugac}):

\begin{itemize}
\item[]{\small
I share the opinion and the same opinion
is shared by my Polish colleagues that the presence of
Aleksandroff, Khinchin, Kolmogorov, and Shnirelman who  confronted their
former teacher in the most dishonest manner and slanderously accused him
is intolerable at any
meeting of decent persons.}
\end{itemize}

The method of political insinuations and slander was
used against the old Muscovite professorship many years before the  {\it Pravda\/} article.
The declaration of November 21, 1930 of the ``initiative group'' of the Moscow Mathematical Society
which consisted of L.~A.~Lyusternik,  L.~G. Shnirelman,
A.~O. Gelfond, and  L.~S. Pontryagin (1908--1988)  and K.~P. Nekrasov
claimed (cp.~\cite{Yush} and \cite{BogRozh}):

\begin{itemize}
\item[]{\small
The acrimonious class struggle in the USSR has pushed the right-wing professoriate into the
camp of the counter-revolution. The reactionary professoriate headed all  sabotaging
organizations and counter-revolutionary parties  that have been disclosed recently.
Owing to the meritorious actions of the {\it OGPU}\endnote{This is the standard abbreviation in Russia of the political police
of those years called the Combined   State Political Department.}
there are divulged the crimes of the whole bunch of scientific bonzes who can
artfully hide themselves under various masks--- from that of cold loyalty to  Soviet
power to highly-advertised profound affection to Soviet power.
 Even in the community of mathematicians some active counter-revolutionaries
 were revealed. Under arrest for participation in a counter-revolutionary organization
 is Professor Egorov, the undisputed leader of the Moscow mathematical school,
 the Chairmen of the Moscow Mathematical Society, the former Director of the Mathematical
 Institute, and the candidate of Moscow mathematics to the Academy of
 Sciences---the very same Egorov, the savior of the academic traditions, whom the proletarian
 student body  had been struggling against for a long time  but whose defence was
 a practically unanimous  decision of the community of Moscow mathematicians.}
\end{itemize}
\noindent

 D.~F. Egorov (1869--1931) was the teacher  of Luzin. Shortly before D.~F. Egorov  had been arrested, and Luzin decided it wise to leave the university
(he was later accused of this removal by his students).

In his life's-description, dated as of the late 1970s, L.~S. Pontryagin
asserted (cp.~\cite[p.~91]{LSP}):

\begin{itemize}
\item[]{\small
The two public actions, in 1936 as regards Luzin and in 1939 as regards
elections, were the important stages of
my uprising as a public person. In my opinion both were the struggle
for  rightful ends.}
\end{itemize}

This is totally inconsistent with the position of Luzin who
wrote in his letter of 1934 to L.~V. Kantorovich (1912--1986) after the
ugly declaration signed by A.~O. Gelfond that his choice in Moscow for the forthcoming election
of corresponding members of the Academy ``will be Gelfond who has recently made
a discovery worthy of a genius'' (cp.~\cite{RK}).

And in 1939 Luzin wrote to V.~I. Vernadsky (cp.~\cite[p.~105]{Antipenko}):

\begin{itemize}
\item[]{\small
Vladimir Ivanovich, the candidates in mathematics---Sobolev and Kolmogorov---are good.
I will vote for them.}
\end{itemize}

A broad campaign against Luzin and ``Luzinism''  waged over this country
in 1936 (cp.~\cite[p.~757--767]{Kol}).
Fortunately, Luzin was not repressed  nor expelled from the Academy.
Some historians opine that there was a relevant oral direction of
I.~V. Stalin.\endnote{It was disclosed recently that the above-mentioned
letter of P.~S. Kapitsa to V.~M. Molotov was multiplied in 16 copies
for the members of the Political Bureau of the All-Union
Communist Party (Bolsheviks) and discussed over with other letters
in support of Luzin.}
But the badge of an enemy under the mask of a Soviet citizen
was pinpointed to Luzin during 14~years up to his death.
The monstrosity over  Luzin is absolutely incomparable with the alleged accusations of
moral misconduct.

\section*{Mathematical Roots of the Luzin Case}

The human passions and follies behind the 1930s tragedy of mathematics in
Russia are obvious: love and hatred, jealousy and admiration, vanity and
modesty, generosity and careerism, etc. But was there a mathematical background?
Some roots are visible.

We are granted the blissful world that has the indisputable property
of unicity. The solitude of reality was perceived by our ancestors as
the ultimate proof of unicity. This argument resided behind the
incessant attempts at proving the fifth postulate of Euclid. The same
gives grounds for the common search of the unique best solution of any
human problem.

Mathematics has never liberated itself from the tethers of
experimentation. The reason is not the simple fact that we still
complete proofs by declaring ``obvious.'' Alive and rather popular are
the views of mathematics as a toolkit for natural sciences. These
stances may be expressed by the slogan ``mathematics is experimental
theoretical physics.'' Not less popular is the dual claim
``theoretical physics is experimental mathematics.'' This short
digression is intended to point to the interconnections of the trains
of thought in mathematics and natural sciences.

It is worth observing that the dogmata of faith and the principles of
theology are also well reflected in the history of mathematical
theories. Variational calculus was invented in search of better
understanding of the principles of mechanics, resting on the religious
views of the universal beauty and harmony of the act of creation.

The twentieth century marked an important twist in the content of
mathematics. Mathematical ideas imbued the humanitarian sphere and,
primarily, politics, sociology, and economics. Social events are
principally volatile and possess a high degree of uncertainty.
Economic processes utilize a wide range of the admissible ways of
production, organization, and management. The nature of nonunicity in
economics transpires: The genuine interests of human beings cannot
fail to be contradictory. The unique solution is an oxymoron in any
nontrivial problem of economics which refers to the distribution of
goods between a few agents. It is not by chance that the social
sciences and instances of humanitarian mentality invoke the numerous
hypotheses of the best organization of production and consumption, the
most just and equitable social structure, the codices of rational
behavior and moral conduct, etc.

The twentieth century became the age of freedom. Plurality and unicity were
confronted as collectivism and individualism. Many particular
phenomena of life and culture reflect their distinction. The
dissolution of monarchism and tyranny was accompanied by the rise of
parliamentarism and democracy. Quantum mechanics and  Heisenberg's
uncertainty incorporated plurality in physics. The waves of modernism
in poetry and artistry should be also listed. Mankind had changed all
valleys of residence and dream.

In mathematics the quest for plurality led to the abandonment of the
overwhelming pressure of unicity and categoricity. The latter ideas were
practically absent, at least minor, in Ancient Greece and sprang to
life in the epoch of absolutism and Christianity. G.~Cantor (1845--1918)
was a~harbinger of mighty changes, claiming that ``das {\it Wesen} der
{\it Mathematik}
liegt gerade in ihrer {\it Freiheit}.'' Paradoxically, the resurrection of
freedom expelled mathematicians from the Cantor paradise.

Nowadays we are accustomed to the unsolvability and undecidability of
many problems. We see only minor difficulties in accepting nonstandard
models and modal logics. It does not worry us that the problem of the
continuum is undecidable within Zermelo--Fraenkel set theory. However
simple nowadays, these stances of thought seemed opportunistic and
controversial at the times of Luzin. The successful breakthroughs of
the great students of Luzin were based on the rejection of his
mathematical ideas. This is a psychological partly Freudian background
of the Luzin case. His gifted students smelled the necessity of
liberation from description and the pertinent blissful dreams of Luzin
which were proved to be undecidable in favor of freedom for
mathematics. His students were misled and consciously or unconsciously
transformed the noble desire for freedom into the primitive hatred and
cruelty. This transformation is a popular fixation and hobby horse
of the human beings through the ages.

Terrible and unbearable is the lightheaded universal fun of putting
the blame entirely on Luzin for the crimes in mathematics  which he was
hardly guilty of with the barely concealed intention to revenge his
genuine and would-be private and personal sins. We should try and
understand that the ideas of description, finitism, intuitionism, and
similar heroic attempts at the turn of the twentieth century in search of
the sole genuine and ultimate foundation were unavoidable by way of
liberating mathematics from the illusionary dreams of categoricity.
The collapse of the eternal unicity and absolutism was a triumph and
tragedy of the mathematical ideas of the first two decades of the last
century. The blossom of the creative ideas of Luzin's students stemmed
partly from his mathematical illusions in description.

The struggle against Luzin had mathematical roots which were
impossible to extract and explicate those days. We see clearly now
that the epoch of probability, functional analysis, distributions,
topology began when the idea of the ultimate unique foundation was
ruined for ever. K.~G\"odel (1906--1978) had explained some trains of thought behind
the phenomenon, but the mathematicians par excellence felt them with
inborn intuition and challenge of mind.

It was Luzin whom the vision of the Moscow schools of today  had started with.
Luzin was interested in foundations, and description for him was the method of
understanding the whole of mathematics. Sprang to life as the theory of measurability,
description has not passed away---it is alive in recursive analysis and other
ideas related to computability ant the Church thesis.

Description plays the same role in regard to finitism and
intuitionism as absolute geometry plays in regard to elliptic and hyperbolic geometries.
The procedures and ideology of description are the forerunners of the ideas of computability
and algorithm. The creative contribution of A.~N. Kolmogorov into algorithm theory,
computability, and complexity comprises the component that stems from description.
Probability theory, turbulence, and analysis constitute the component that is rooted in
the refusal from description. Mathematics reduces to neither finitism, nor
intuitionism, nor description. It is not categorical---it is free.
In the twentieth century the freedom of mathematics was best demonstrated in Russia
by A.~N. Kolmogorov, a student of Luzin and the teacher of new generations of mathematicians
in this country.  He harbored more freedom in mathematics, which made him a greater mathematician.

It is the tragedy of  mathematics in Russia  that the  noble
endeavor for  freedom  had launched the political
monstrosity of the scientific giants disguised into the
cassocks of Torquemada.

\section*{A Few Lessons}

History and decedents are out of the courts of justice.
Scientists and ordinary persons must see  and collect facts.
Never accuse the passed away, but calmly and openly point out
that which was in reality. Explain  the difference between
moral accusations and political insinuations to the youth.
Demonstrate the difficulty and necessity  of the repairing
of mistakes  and repentance. Show how easy it is
to forgive oneself and accuse the others.

We must work out and transfer to the next generations
the objective views of the past. Of its successes and tragedies.
With love and doubts, with the understanding of our unfortunate
fate and the honor of objectivity. It is the personal faults and
failures that we are to accuse and repair first of all.
They knew even in Ancient Rome that we should  tell nothing or good about the dead.
Facts did never pass away. Luzin was accused by the community
of Moscow mathematicians as well as by the Academy of Sciences.

The defence of Stalinism consists often in proclaiming that
the ugly misdeeds and crimes of Stalinism are the fault of Stalin,
Beriya, Mekhlis, Kolman, and the hord of their minor clones
and replicas, while in fact Stalinism was created by millions.
Stalinism in science was mainly raised by scientists themselves.
It is indecent to pretend that Luzin's students safeguarded their teacher and
science from Stalinism. Luzin was a victim of social ostracism
and lived with the brand and stigma of an adversary with a Soviet mask during fourteen
years up to his death. He became an exemplary outcast for  Stalinism---an adversary
at large. Luzin's colleagues and students lowered and neglected him, which culminated in
slapping in his face and spitting on  his tomb.
Luzin has passed away but the false accusations in sabotage and obsequiousness
are still  effective.

We cannot forget the words of Luzin (cp.~\cite[p.~73]{Case}):

\begin{itemize}
\item[]{\small
As regards the last paragraph of the article in {\it Pravda\/} where
the monstrous accusations were made against me in servitude to
the present-day masters
of faschistoid science, I state with  the full understanding of my political
responsibility of a scientist of  international reputation and a citizen of
the USSR that the Editorial Board of {\it Pravda\/} was deliberately
misled to delusion by the persons who had informed it about this.
This is refuted by my whole life and activities as a scientist and  a person.}
\end{itemize}

Any attempt at discerning  morality in the past immorality
is dangerous since it feeds this immorality by creating
the comfortable environment of immorality in the present and future.
The stamina of a scientist by belief is a discontinuous function.
Science does not inoculate morality. 
Evil and genius coexist from time to time. 

History is a branch of science, but  science goes hand in hand with conscience.
Therefore history  is not only a scientific discipline 
but also a matter of responsibility of the present.
No one can change the history but we, the humankind,  are not 
indifferent  onlookers  on the past.  History is not anything existant without people.
The past is the past of the present and as such it is part 
of the present-day responsibility.  The alive rather than the dead  
are responsible for what was done in the past. It is we who
change life so creating the future. History calls upon our conscience and 
hits in our hearts alike the ashes of Claes' did in Till Eulenspiegel's.

\section*{Acknowledgements}

This essay summarizes the author's position at the informal lobby discussions at
the Mini-Symposium on Convex Analysis which was held at Lomonosov Moscow
State University, February 2-4 (2007). The author is especially grateful to Professor
V.~M. Tikhomirov, Chairman of the Mini-Symposium, for endurance, friendliness, and hospitality.

This article is an expanded and updated version of the author's article
``Roots of Luzin's case.''\endnote{J.~Appl.~Indust. Math., {\bf 1}: 3, 261--267 (2007).}
The author have received many comments and replies with criticism and advice  aiming at
 improvement of the article, as well as indication of many new documentary sources and
reminiscences. The author thanks all his readers and correspondents.

\bibliographystyle{plain}

\theendnotes

\end{document}